\documentclass[leqno]{article}
\usepackage{amsmath}
\usepackage{amssymb}
\usepackage{dsfont}
\usepackage[usenames]{color}
\usepackage[colorlinks]{hyperref}
\usepackage{maple2e}
\DefineParaStyle{Maple Output}
\DefineParaStyle{Maple Plot}
\DefineParaStyle{Warning}
\DefineCharStyle{2D Math}
\DefineCharStyle{2D Output}

\def\Reg{\operatorname{Reg}}

\begin{document}

\title{\sc Regular integers modulo $n$}
\author{{\bf L\'aszl\'o T\'oth} (P\'ecs, Hungary)}
\date{}
\maketitle

\centerline{Annales Univ. Sci. Budapest., Sect. Comp., {\bf 29}
(2008), 263-275} \vskip4mm

\vskip1mm
\centerline{\it Dedicated to Professor Imre K\'atai on his 70th birthday}
\vskip1mm

\begin{abstract} Let $n=p_1^{\nu_1}\cdots p_r^{\nu_r} >1$ be an integer.
An integer $a$ is called regular (mod $n$) if there is an integer
$x$ such that $a^2x\equiv a$ (mod $n$). Let $\varrho(n)$ denote the
number of regular integers $a$ (mod $n$) such that $1\le a\le n$.
Here $\varrho(n)=(\phi(p_1^{\nu_1})+1)\cdots
(\phi(p_r^{\nu_r})+1)$, where $\phi(n)$ is the Euler function.  In
this paper we first summarize some basic properties of regular
integers (mod $n$). Then in order to compare the rates of growth of
the functions $\varrho(n)$ and $\phi(n)$ we investigate the average
orders and the extremal orders of the functions
$\varrho(n)/\phi(n)$, $\phi(n)/\varrho(n)$ and $1/\varrho(n)$.
\end{abstract}



\vskip1mm {\bf 1. Introduction}

\vskip1mm Let $n>1$ be an integer. Consider the integers $a$ for
which there exists an integer $x$ such that $a^2x\equiv a$ (mod
$n$). In the background of this property is that an element $a$ of a
ring $R$ is said to be regular (following J. von Neumann) if there
is an $x\in R$ such that $a=axa$. In case of the ring $\mathds{Z}_n$
this is exactly the condition of above.

Properties of these integers were investigated by J. Morgado
\cite{M72},  \cite{M74}, who called them regular (mod $n$). In a
recent paper O. Alkam and E. A. Osba \cite{AO} using ring theoretic
considerations rediscovered some of the statements proved
elementarly by J. Morgado. It was observed in \cite{M72}, \cite{M74}
that $a>1$ is regular (mod $n$) if and only if the gcd $(a,n)$ is a
unitary divisor of $n$. We recall that $d$ is said to be a unitary
divisor of $n$ if $d\mid n$ and gcd $(d,n/d)=1$, notation $d \mid
\mid n$.

These integers occur in the literature also in an other context. It
is said that an integer $a$ possesses a weak order (mod $n$) if
there exists an integer $k\ge 1$ such that $a^{k+1} \equiv a$ (mod
$n$). Then the weak order of $a$ is the smallest $k$ with this
property, see \cite{J81}, \cite{F}. It turns out that $a$ is regular
(mod $n$) if and only if $a$ possesses a weak order (mod $n$).

Let $\Reg_n=\{a: 1\le a\le n$, $a$ is regular (mod $n$)$\}$ and let
$\varrho(n)=\# \Reg_n$ denote the number of regular integers $a$
(mod $n$) such that $1\le a\le n$. This function is multiplicative
and $\varrho(p^{\nu})=\phi(p^{\nu})+1= p^{\nu}-p^{\nu-1}+1$ for
every prime power $p^{\nu}$ ($\nu \ge 1$), where $\phi$ is the Euler
function. Consequently, $\displaystyle \varrho(n)=\sum_{d \mid \mid
n} \phi(d)$ for every $n\ge 1$, also $\phi(n)< \varrho(n)\le n$ for
every $n>1$, and $\varrho(n)= n$ if and only if $n$ is squarefree,
see \cite{M72}, \cite{J81}, \cite{AO}.

Let us compare the functions $\varrho(n)$ and $\phi(n)$. The first
few values of $\varrho(n)$ and $\phi(n)$ are given by the next
tables ($\varrho(n)$ is sequence $A055653$  in Sloane's On-Line
Encyclopedia of Integer Sequences \cite{S}). Note that $\varrho(n)$
is even iff $n\equiv 2$ (mod $4$), and $\sqrt{n}\le \varrho(n)\le n$
for every $n\ge 1$, see \cite{AO}.

$$\vbox{\offinterlineskip
\halign{\strut \quad
      \hfil $\ # \ $ \hfil &\vrule \vrule
      \hfil $\ # \ $ \hfil &\vrule    \hfil $\ #\ $ \hfil &\vrule
      \hfil $\ # \ $ \hfil &\vrule    \hfil $\ #\ $ \hfil &\vrule
      \hfil $\ # \ $ \hfil &\vrule    \hfil $\ #\ $ \hfil &\vrule
      \hfil $\ # \ $ \hfil &\vrule    \hfil $\ #\ $ \hfil &\vrule
      \hfil $\ # \ $ \hfil &\vrule    \hfil $\ #\ $ \hfil &\vrule
      \hfil $\ # \ $ \hfil &\vrule    \hfil $\ #\ $ \hfil &\vrule
      \hfil $\ # \ $ \hfil &\vrule    \hfil $\ #\ $ \hfil &\vrule
      \hfil $\ # \ $ \hfil    \cr
      n & 1 & 2 & 3 & 4 & 5 & 6 & 7 & 8 & 9 &
               10 & 11 & 12 & 13 & 14 & 15  \cr
      \noalign{\hrule} \noalign{\hrule}
      \varrho(n) & 1 & 2 & 3 & 3 & 5 & 6 & 7 & 5 &
                7 & 10 & 11 & 9 & 13 & 14 & 15  \cr \noalign{\hrule} \noalign{\hrule}
      \phi(n) & 1 & 1 & 2 & 2 & 4 & 2 & 6 & 4 &
                6 & 4 & 10 & 4 & 12 & 6 & 8
\cr}}
$$

$$\vbox{\offinterlineskip
\halign{\strut \quad
      \hfil $\ # \ $ \hfil &\vrule \vrule
      \hfil $\ # \ $ \hfil &\vrule    \hfil $\ #\ $ \hfil &\vrule
      \hfil $\ # \ $ \hfil &\vrule    \hfil $\ #\ $ \hfil &\vrule
      \hfil $\ # \ $ \hfil &\vrule    \hfil $\ #\ $ \hfil &\vrule
      \hfil $\ # \ $ \hfil &\vrule    \hfil $\ #\ $ \hfil &\vrule
      \hfil $\ # \ $ \hfil &\vrule    \hfil $\ #\ $ \hfil &\vrule
      \hfil $\ # \ $ \hfil &\vrule    \hfil $\ #\ $ \hfil &\vrule
      \hfil $\ # \ $ \hfil &\vrule    \hfil $\ #\ $ \hfil &\vrule
      \hfil $\ # \ $ \hfil    \cr
      n & 16 & 17 & 18 & 19 & 20 & 21 & 22 & 23 & 24 &
               25 & 26 & 27 & 28 & 29 & 30  \cr
      \noalign{\hrule} \noalign{\hrule}
      \varrho(n) & 9 & 17 & 14 & 19 & 15 & 21 & 22 & 23 &
                15 & 21 & 26 & 19 & 21 & 29 & 30   \cr \noalign{\hrule} \noalign{\hrule}
      \phi(n) & 8 & 16 & 6 & 18 & 8 & 12 & 10 & 22 &
                8 & 20 & 12 & 18 & 12 & 28 & 8
 \cr}}
$$

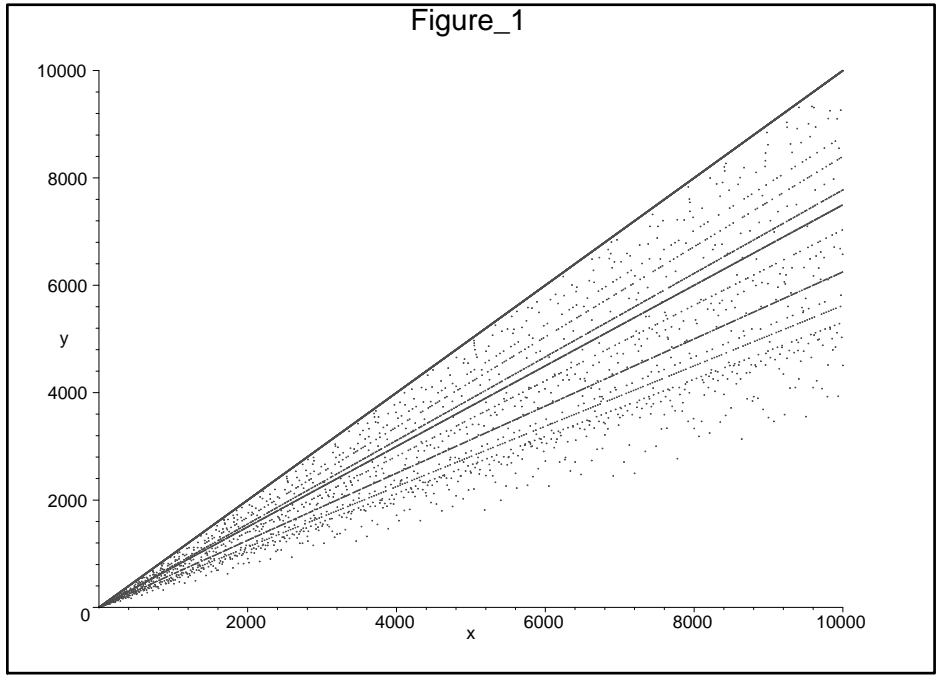
Figure 1 is a plot of the function $\varrho(n)$ for $1\le n\le 10\, 000$.

{\mapleplot{Fig101.eps}}

For the Euler $\phi$-function,
\[
\lim_{x\to \infty} \frac1{x^2} \sum_{n\le x} \phi(n) =
\frac{3}{\pi^2} \approx 0.3039.
\]

The average order of the function $\varrho(n)$ was considered in
\cite{J81}, \cite{F}. One has
\[
\lim_{x\to \infty} \frac1{x^2} \sum_{n\le x} \varrho(n) =\frac1{2} A
 \approx 0.4407,
\]
where \[ A=\prod_p \left(1-\frac1{p^2(p+1)}\right)=\zeta(2)\prod_p
\left(1-\frac1{p^2}-\frac1{p^3}+\frac1{p^4}\right)\approx 0.8815
\]
is the so called quadratic class-number constant. For its numerical
evaluation see \cite{NM}.

More exactly,
\[
\sum_{n\le x} \varrho(n)=\frac1{2}Ax^2+ R(x),
\]
where $R(x)=O(x\log^3 x)$, given in \cite{J81} using elementary
arguments. This was improved into $R(x)=O(x\log^2 x)$ in \cite{Y86},
and into $R(x)=O(x\log x)$ in \cite{HS92}, using analytic methods.
Also, $R(x)=\Omega_{\pm}(x\sqrt{\log \log x})$, see \cite{HS92}.

In this paper we first summarize some basic properties of regular
integers (mod $n$). We give also their direct proofs, because the
proofs of \cite{M72}, \cite{M74} are lengthy and those of \cite{AO}
are ring theoretical.

Then in order to compare the rates of growth of the functions
$\varrho(n)$ and $\phi(n)$ we investigate the average orders and the
extremal orders of the functions $\varrho(n)/\phi(n)$,
$\phi(n)/\varrho(n)$ and $1/\varrho(n)$. The study of the minimal
order of $\varrho(n)$ was initiated in \cite{AO}.

\vskip1mm {\bf 2. Characterization of regular integers (mod $n$)}

\vskip1mm The integer $a=0$ and those coprime to $n$ are regular (mod $n$) for
each $n>1$. If $a\equiv b$ (mod $n$), then $a$ and $b$ are regular
(mod $n$) simultaneously. If $a$ and $b$ are regular (mod $n$),
then $ab$ is also regular (mod $n$).

In what follows let $n>1$ be of canonical form $n=p_1^{\nu_1}\cdots
p_r^{\nu_r}$.

\vskip1mm {\bf Theorem 1.} {\it For an integer $a\ge 1$ the following
assertions are equivalent:

i) $a$ is regular (mod $n$),

ii) for every $i\in \{1,...,r\}$ either $p_i\nmid a$ or
$p_i^{\nu_i}\mid a$,

iii) $(a,n)=(a^2,n)$,

iv) $(a,n) \mid \mid n$,

v) $\displaystyle a^{\phi(n)+1}\equiv a$ (mod $n$),

vi) there exists an integer $k\ge 1$ such that $a^{k+1} \equiv a$
(mod $n$).}

\vskip1mm {\bf Proof.} i) $\Rightarrow$ ii). If $a^2x\equiv a$ (mod
$n$) for an integer $x$, then $a(ax-1)\equiv 0$ (mod $p_i^{\nu_i}$)
for every $i$. We have two cases: $p_i\nmid a$ and $p_i\mid a$. In
the second case, since $(a,ax-1)=1$, obtain that $a \equiv 0$ (mod
$p_i^{\nu_i}$).

ii) $\Rightarrow$ i). If $p_i^{\nu_i}\mid a$, then $a^2x \equiv a$
(mod $p_i^{\nu_i}$) for any $x$. If $p_i\nmid a$, then the linear
congruence $ax\equiv 1$ (mod $p_i^{\nu_i})$ has solutions in $x$ and
obtain also $a^2x \equiv a$ (mod $p_i^{\nu_i}$).

ii) $\Leftrightarrow$ iii). Follows at once by the property of the
gcd.

ii) $\Leftrightarrow$ iv) Follows at once by the definition of the
unitary divisors (the unitary divisors of a prime power $p^\nu$ are
$1$ and $p^\nu$).

ii) $\Rightarrow$ v) (\cite{AO}) If $p_i^{\nu_i}\mid a$, then
$a^{\phi(n)+1} \equiv a$ (mod $p_i^{\nu_i}$). If $p_i\nmid a$, then
using Euler's theorem, $a^{\phi(n)+1}\equiv
a(a^{\phi(p_i^{\nu_i})})^{\phi(n)/\phi(p_i^{\nu_i})} \equiv a$ (mod
$p_i^{\nu_i})$. Therefore $a^{\phi(n)+1}\equiv a$ (mod
$p_i^{\nu_i}$) for every $i$ and $a^{\phi(n)+1}\equiv a$ (mod $n$).

v) $\Rightarrow$ i) (\cite{AO}) If $a^{\phi(n)+1}\equiv a$ (mod
$n$), then $a^2 a^{\phi(n)-1}\equiv a$ (mod $n$), hence $a^2x\equiv
a$ (mod $n$) is verified for $x=a^{\phi(n)-1}$ (which is the von
Neumann inverse of $a$ in $\mathds Z_n$).

v) $\Rightarrow$ vi) Immediate by taking $k=\phi(n)$.

vi) $\Rightarrow$ i) If $a^{k+1}\equiv a$ (mod $n$) for an integer $k\ge 1$,
then $a^2x\equiv a$ (mod $n$) holds for $x=a^{k-1}$, finishing the
proof.

\vskip1mm Note that the proof of i) $\Leftrightarrow$ v) given in
\cite{M74} uses Dirichlet's theorem on arithmetic progressions,
which is unnecessary.

\vskip1mm {\bf Theorem 2.} {\it The function $\varrho(n)$ is
multiplicative and $\varrho (p^\nu)=p^\nu-p^{\nu-1}+1$ for every
prime power $p^\nu$ ($\nu \ge 1$). For every $n\ge 1$,
\[
\varrho(n)=\sum_{d \mid \mid n} \phi(d).
\]}

\vskip1mm {\bf Proof.} By Theorem 1, $a$ is regular (mod $n$) iff
for every $i\in \{1,...,r\}$ either $p_i\nmid a$ or $p_i^{\nu_i}\mid
a$.

Let $a\in \Reg_n$. If $p_i\nmid a$ for every $i$, then $(a,n)=1$,
the number of these integers $a$ is $\phi(n)$. Suppose that
$p_i^{\nu_i}\mid a$ for exactly one value $i$ and that for all $j\ne
i$, $(p_j,a)=1$. Then $a=bp_i^{\nu_i}$, where $1\le b\le
n/p_i^{\nu_i}$ and $(b, n/p_i^{\nu_i})=1$. The number of such
integers $a$ is $\phi(n/p_i^{\nu_i})$. Now suppose that
$p_i^{\nu_i}\mid a$, $p_j^{\nu_j}\mid a$, $i<j$, and for all $k\ne
i, k\ne j$, $(p_i,a)=(p_j,a)=1$. Then $a=cp_i^{\nu_i}p_j^{\nu_j}$,
where $1\le c\le n/(p_i^{\nu_i}p_j^{\nu_j})$ and $(c,
n/(p_i^{\nu_i}p_j^{\nu_j}))=1$. The number of such integers $a$ is
$\phi(n/(p_i^{\nu_i}p_j^{\nu_j}))$, etc. We obtain
\[
\varrho(n) = \phi(n)+\sum_{1\le i\le r} \phi(n/p_i^{\nu_i})+ \sum_{1\le i<j\le r}
\phi(n/p_i^{\nu_i}p_j^{\nu_j})+...+
\phi(n/(p_1^{\nu_1}\cdots p_r^{\nu_r})). \]

Let $y_i=\phi(p_i^{\nu_i})$, $1\le i\le r$, and $y=y_1\cdots y_r$. Then $\phi(n)=y$ and
\[
\varrho(n) = y+ \sum_{1\le i\le r} \frac{y}{y_i} + \sum_{1\le i<j\le r}
\frac{y}{y_iy_j}+...+
\frac{y}{y_1\cdots y_r} = \]
\[ =(y_1+1)\cdots (y_r+1)= (\phi(p_1^{\nu_i})+1)\cdots (\phi(p_r^{\nu_r})+1). \]

The given representation of $\varrho(n)$ now follows at once taking into account that the
unitary convolution
preserves the  multiplicativity of functions, see for example \cite{Mc86}.

Another method, see \cite{M72}: Group the integers $a\in \{
1,2,...,n\}$ according to the value $(a,n)$. Here $(a,n)=d$ if and
only if $(j,n/d)=1$, where $a=jd$, $1\le j\le n/d$, hence the number
of integers $a$ with $(a,n)=d$ is $\phi(n/d)$. According to Theorem
1, $a$ is regular (mod $n$) if and only if $d=(a,n)\mid \mid n$, and
obtain that
\[
\varrho(n)=\sum_{d \mid \mid n} \phi(n/d)= \sum_{d \mid \mid n} \phi(d).
\]

Now the multiplicativity of $\varrho(n)$ is a direct consequence of this representation.

\vskip1mm Let $S(n)$ denote the sum of regular integers $a\in
\Reg_n$. We give a simple formula for $S(n)$, not considered in the
cited papers, which is analogous to $\sum_{1\le a\le n, (a,n)=1} a
=n\phi(n)/2$ ($n>1$).

\vskip1mm {\bf Theorem 3.} {\it For every $n\ge 1$,
\[
S(n)=\frac{n(\varrho(n)+1)}{2}.
\]}

\vskip1mm {\bf Proof.} Similar to the counting procedure of above or
by grouping the integers $a\in \{ 1,2,...,n\}$ according to the
value $(a,n)$:
\[
S(n)= \sum_{a\in \Reg_n} a =
\sum_{d\mid \mid n} \sum_{\substack{a\in \Reg_n \\ (a,n)=d}} a = \sum_{d\mid \mid n} d
\sum_{\substack{j=1\\(j,n/d)=1}}^{n/d} j=
\]
\[
=n+ \sum_{\substack{d\mid \mid n\\ d<n}} d\frac{n\phi(n/d)}{2d} = n+
\frac{n}{2} \sum_{\substack{d\mid \mid n\\ d<n}} \phi(n/d)= \frac{n(\varrho(n)+1)}{2}.
\]

\vskip1mm {\bf 3. Average orders}

\vskip1mm {\bf Theorem 4.} {\it For the quotient
$\varrho(n)/\phi(n)$ we have
\[ \sum_{n\le x} \frac{\varrho(n)}{\phi(n)} = B x + O(\log^2 x),\]
where $B=\pi^2/6\approx 1.6449$. }

\vskip1mm {\bf Proof.} By Theorem 2, $\varrho(p^{\nu})/\phi(p^{\nu})
= 1+1/\phi(p^{\nu})$ for every prime power $p^{\nu}$ ($\nu \ge 1$).
Hence, taking into account the multiplicativity, for every $n\ge 1$,
\[
\frac{\varrho(n)}{\phi(n)} =\sum_{d\mid \mid n} \frac1{\phi(d)}.
\]

Using this representation (given also in \cite{AO}) we obtain
\[
\sum_{n\le x} \frac{\varrho(n)}{\phi(n)} =\sum_{\substack{de\le x\\
(d,e)=1}} \frac1{\phi(d)} = \sum_{d\le x}
\frac1{\phi(d)}\sum_{\substack{e\le x/d \\ (e,d)=1}} 1= \]
\[
=\sum_{d\le x} \frac1{\phi(d)}
\left(\frac{\phi(d)x}{d^2}+O(2^{\omega(d)}) \right)= x \sum_{d\le x}
\frac1{d^2} + O\left(\sum_{d\le x} \frac{2^{\omega(d)}}{\phi(d)}
\right),
\]
where $\omega(d)$ denotes, as usual, the number of distinct prime
factors of $d$. Furthermore, let $\tau(n)$ and $\sigma(n)$ denote
the number and the sum of divisors of $n$, respectively. Using that
$\phi(n)\sigma(n)\gg n^2$, we have $2^{\omega(d)}/\phi(d) \ll
\tau(d)\sigma(d)/d^2$. Here $\sum_{d\le x} \tau(d)\sigma(d) \ll
x^2\log x$, according to a result of Ramanujan, and obtain by
partial summation that the error term is $O(\log^2 x)$.

Figure 2 is a plot of the error term $\sum_{n\le x} \varrho(n)/\phi(n) - B x$ for $1\le x
\le 1000$.

{\mapleplot{Fig201.eps}}

Consider now the quotient $f(n)=\phi(n)/\varrho(n)$, where
$f(n)\le 1$. According to a well-known result of H. Delange, $f(n)$
has a mean value given by
\[
C= \prod_p \left(1-\frac1{p}\right)\left(1+ \sum_{\nu=1}^{\infty}
\frac{f(p^{\nu})}{p^{\nu}} \right)= \prod_p \left(1-\frac1{p}\right)
\left(1+ (1-\frac1{p}) \sum_{\nu=1}^{\infty}
\frac1{p^{\nu}-p^{\nu-1}+1} \right).
\]

Here $C\approx 0.6875$, which can be obtained using that for every
$k\ge 1$,
\[ C= \prod_p \left(1-\frac1{p}\right)\left(1+ (1-\frac1{p})
\sum_{\nu=1}^k \frac1{p^{\nu}-p^{\nu-1}+1} +\frac1{p^kr_p}\right),
\]
where $p-1<r_p<p$ for each prime $p$.

We prove the following asymptotic formula:

\vskip1mm {\bf Theorem 5.} {\it
\[
\sum_{n\le x} \frac{\phi(n)}{\varrho(n)}= C x + O((\log x)^{5/3}(\log \log x)^{4/3}).
\]}

\vskip1mm {\bf Proof.} For $f(n)=\phi(n)/\varrho(n)$ let
\[
f(n) =\sum_{d\mid n} \frac{\phi(d)}{d}\, v(n/d),
\]
that is, in terms of the Dirichlet convolution, $f=\phi/E * v$,
$f=\mu/E*I*v$, $v=f* I/E * \mu$, where $\mu(n)$ is the M\"obius
function, $E(n)=n$, $I(n)=1$ ($n\ge 1$).

The function $v(n)$ is multiplicative, for every prime power $p^\nu$
($\nu \ge 1$),
\[
v(p^{\nu})=f(p^{\nu})-\left(1-\frac1{p}\right) \left(f(p^{\nu-1})+
\frac1{p}f(p^{\nu-2})+... +\frac1{p^{\nu-2}}f(p)+
\frac1{p^{\nu-1}}\right),
\]
and $v(p)=0$, $|v(p^2)|\le 1/p$ for every prime $p$. Also,
\[
f(p^{\nu})=\frac{p^{\nu}-p^{\nu-1}}{p^{\nu}-p^{\nu-1}+1}=\frac{1-1/p}{1-(1/p-
1/p^{\nu})}=
\left(1-\frac1{p}\right)\left(1+\left(\frac1{p}-\frac1{p^\nu}\right)+
\left(\frac1{p}-\frac1{p^\nu}\right)^2+...\right)=
\]
\[=
\left(1-\frac1{p}\right)\left(1+\frac1{p}+\frac1{p^2}
+...+\frac1{p^{\nu}}- \frac1{p^{\nu}} +O\left(\frac1{p^{\nu+1}}\right)\right),
\]
and obtain that for every fixed $\nu \ge 3$,
\[
f(p^{\nu})= 1-\frac1{p^{\nu}}+O\left(\frac1{p^{\nu+1}}\right),
\]
consequently,
\[
v(p^\nu)= 1-\frac1{p^\nu}+O\left(\frac1{p^{\nu+1}}\right)-\left(1-\frac1{p}\right)\left(1+\frac1{p}+...+\frac1{p^{\nu-1}} - \frac{\nu-1}{p^{\nu-1}} +O\left(\frac1{p^{\nu}}\right)\right),
\]
\[
v(p^\nu)= \frac{\nu-1}{p^{\nu-1}} + O\left(\frac1{p^\nu}\right). \leqno(*)
\]

It follows that there exists $x_0$ such that for every prime $p>x_0$ and for every
$\nu \ge 3$,
\[
|v(p^\nu)|\le \frac1{p^{3\nu/5}}. \leqno(**)
\]

Now we show that
\[
\sum_{n\le x} v(n)=O(\log x), \quad \sum_{n>x} \frac{v(n)}{n} =
O\left(\frac{\log x}{x}\right).
\]

We deduce the first estimate, the second one will follow by partial
summation. Let ${\cal M}_1 =\{ n: \ p \mid n \ \Rightarrow \ p\le
x_0\}$, ${\cal M}_2 =\{ n: \ p \mid n \ \Rightarrow \ p^3\mid n,
p>x_0\}$, ${\cal M}_3 =\{ n: \ p \mid n \ \Rightarrow \ p^2\mid n,
p^3\nmid n, p>x_0\}$. If $v(n)\ne 0$, then $n$ can be written
uniquely as $n=n_1n_2n_3$, where $n_1\in {\cal M}_1$, $n_2\in {\cal
M}_2$, $n_3\in {\cal M}_3$. We have the following estimates.

If $n_3\in {\cal M}_3$, then $n_3=m^2$ with $|\mu(m)|=1$. Using
$|v(p^2)|\le 1/p$ we have $|v(n_3)|\le 1/m$, and
\[
\sum_{\substack{n_3\le x\\ n_3\in {\cal M}_3}} v(n_3) \ll \sum_{m\le
\sqrt{x}} \frac{|\mu(m)|}{m}\ll \log x.
\]

By (**), for $x_0$ sufficiently large,
\[ \sum_{\substack{n_2\le x\\ n_2\in {\cal M}_2}} v(n_2) \ll
\prod_{p>x_0} \left(1+|v(p^3)|+|v(p^4)|+...\right)\ll \prod_{p>
x_0} \left(1+\frac1{p^{9/5}}+\frac1{p^{12/5}}+ ...\right) \ll \] \[
\ll \prod_{p> x_0} \left(1+\frac{2}{p^{9/5}}\right)< \infty.
\]

Using (*) we also have
\[ \sum_{\substack{n_1\le x\\ n_1\in {\cal M}_1}} v(n_1) \ll
\prod_{p\le x_0} \left(1+\frac1{p}+ |v(p^3)|+ |v(p^4)|+...\right)<
\infty.
\]

Hence
\[
\sum_{n\le x} v(n)= \sum_{n_1n_2n_3\le x} |v(n_1)|\, |v(n_2)|\,
|v(n_3)| = \sum_{n_1n_2 \le x} |v(n_1)|\, |v(n_2)|\, \sum_{n_3\le
x/n_1n_2} |v(n_3)| \ll \log x. \]

Now applying the following well-known result of Walfisz,
\[
\sum_{n\le x} \frac{\phi(n)}{n} =\frac{6}{\pi^2} x+ O((\log
x)^{2/3}(\log \log x)^{4/3})
\]
we have
\[
\sum_{n\le x} f(n)= \sum_{d\le x} v(d) \sum_{e\le x/d}
\frac{\phi(e)}{e} = \frac{6}{\pi^2} x \sum_{d\le x} \frac{v(d)}{d} +
O((\log x)^{2/3}(\log \log x)^{4/3}\sum_{d\le x} v(d))= \] \[
=\frac{6}{\pi^2} x \sum_{d=1}^{\infty} \frac{v(d)}{d} + O((\log
x)^{5/3}(\log \log x)^{4/3}),
\]
ending the proof of Theorem 5.

Figure 3 is a plot of the error term $\sum_{n\le x} \phi(n)/\varrho(n) - C x$ for $1\le x\le 1000$.

{\mapleplot{Fig301.eps}}

{\bf Theorem 6.} {\it
\[
\sum_{n\le x} \frac1{\varrho(n)}=D\log x + E
+O\left(\frac{\log^9x}{x}\right),
\]
where $D$ and $E$ are constants,
\[
D=\frac{\zeta(2)\zeta(3)}{\zeta(6)} \prod_p
\left(1-\frac{p(p-1)}{p^2-p+1}\sum_{\nu=1}^{\infty}
\frac1{p^{\nu}(p^{\nu}-p^{\nu-1}+1)} \right).\] }

\vskip1mm {\bf Proof.} Write
\[
\frac1{\varrho(n)}=\sum_{\substack{de=n\\(d,e)=1}}
\frac{h(d)}{\phi(e)},\] where $h$ is multiplicative and for every
prime power $p^{\nu}$ ($\nu \ge 1$),
\[
\frac1{\varrho(p^{\nu})}= h(p^{\nu})+ \frac1{\phi(p^{\nu})}, \quad
h(p^{\nu})=- \frac1{\phi(p^{\nu})(\phi(p^{\nu})+1)},
\]
therefore $h(n)\ll 1/\phi^2(n)$. We need the following known result,
cf. for example \cite{MV07}, p. 43,
\[
\sum_{\substack{n\le x\\ (n,k)=1}} \frac1{\phi(n)} = K a(k)
\left(\log x + \gamma + b(k)\right)+ O\left(2^{\omega(k)} \frac{\log
x}{x}\right),
\]
where $\gamma$ is Euler's constant,
\[
K= \frac{\zeta(2)\zeta(3)}{\zeta(6)}, \ a(k)=\prod_{p\mid
k}\left(1-\frac{p}{p^2-p+1}\right) \le \frac{\phi(k)}{k}, \]
\[
b(k)=\sum_{p\mid k} \frac{\log p}{p-1}- \sum_{p\nmid k} \frac{\log
p}{p^2-p+1} \ll \frac{\psi(k)\log k}{\phi(k)}, \ \text{ with } \
\psi(k)= k \prod_{p\mid k} \left(1+\frac1{p}\right).
\]

We have
\[
\sum_{n\le x} \frac1{\varrho(n)}=\sum_{d\le x} h(d)
\sum_{\substack{e\le x/d\\ (e,d)=1}} \frac1{\phi(e)}=
\]
\[
= K \left((\log x+\gamma)\sum_{d\le x} h(d)a(d)+\sum_{d\le x}
h(d)a(d)(b(d)-\log d) \right)+O\left(\frac{\log x}{x}\sum_{d\le x} d
|h(d)|2^{\omega(d)} \right),
\]
and we obtain the given result with the constants
\[
D=K\sum_{n=1}^{\infty} h(n)a(n), \ E=K\gamma \sum_{n=1}^{\infty}
h(n)a(n) + K \sum_{n=1}^{\infty} h(n)a(n)(b(n)-\log n),
\]
these series being convergent taking into account the estimates of
above. For the error terms,
\[
\sum_{n>x} |h(n)|a(n) \ll \sum_{n>x} \frac1{n\phi(n)} \ll \sum_{n>x}
\frac{\sigma(n)}{n^3}\ll \frac1{x}, \ \sum_{n>x} |h(n)| a(n)\log n
\ll \frac{\log x}{x},
\]
\[
\sum_{n>x} |h(n)a(n)b(n)| \ll \sum_{n>x} \frac{\tau^3(n)\log n}{n^2} \ll
\frac{\log^8 x}{x},
\]
using that $\sum_{n\le x} \tau^3(n) \ll x \log^7 x$ (Ramanujan), and
\[
\sum_{n\le x} n|h(n)|2^{\omega(n)} \ll \sum_{n\le x}
\frac{\tau^3(n)}{n} \ll \log^8 x.
\]

\vskip1mm {\bf 4. Extremal orders}

\vskip1mm Since $\varrho(n)\le n$ for every $n\ge 1$ and $\varrho(p)=p$ for
every prime $p$, it is immediate that \\ $\limsup_{n\to \infty}
\varrho(n)/n = 1$. The minimal order of $\varrho(n)$ is also the same
as that of $\phi(n)$, namely,

\vskip1mm {\bf Theorem 7.} {\it
\[\liminf_{n\to \infty}
\frac{\varrho(n)\log \log n}{n}=e^{-\gamma}.
\]}

\vskip1mm {\bf Proof.} We apply the following result (\cite{TW},
Corollary 1): If $f$ is a nonnegative real-valued multiplicative
arithmetic function such that for each prime $p$,

i) $\rho(p):=\sup_{\nu \ge 0} f(p^{\nu})\le (1-1/p)^{-1}$, and

ii) there is an exponent $e_p=p^{o(1)}\in \mathds N$ satisfying
$f(p^{e_p})\ge 1+1/p$,

then
\[
\displaystyle \limsup_{n\to \infty} \frac{f(n)}{\log \log
n}=e^{\gamma}\prod_p \left(1-\frac1{p}\right) \rho(p).
\]

Take $f(n)=n/\varrho (n)$, where $f(p^{\nu})=(1-1/p+ 1/p^{\nu})^{-1}
< (1-1/p)^{-1}=\rho(p)$, and for $e_p=3$,
\[
f(p^3)> 1+\frac{p^2-1}{p^3-p^2+1}> 1+\frac1{p} \] for every prime
$p$.

\vskip1mm It is immediate that $\liminf_{n\to \infty}
\varrho(n)/\phi(n)=1$. The maximal order of $\varrho(n)/\phi(n)$ is
given by

\vskip1mm {\bf Theorem 8.} {\it
\[
\limsup_{n\to \infty} \frac{\varrho(n)}{\phi(n)\log \log
n}=e^{\gamma}.\] }

\vskip1mm {\bf Proof.} Now let $f(n)=\varrho (n)/\phi(n)$ in the
result given above. Here
\[
f(p^{\nu})=1+\frac1{p^\nu-p^{\nu-1}}\le 1+\frac1{p-1}=
\left(1-\frac1{p}\right)^{-1}= \rho(p), \]
and for $e_p=1$, $f(p)>1+1/(p-1)> 1+1/p$ for every prime $p$.

\vskip1mm {\bf 5.} The plots were produced using Maple. The function $\varrho(n)$ was
generated by the following procedure:
\begin{verbatim}
   rho:= proc(n) local x, i: x:= 1: for i from 1 to nops(ifactors(n)[ 2 ]) do
   p_i:=ifactors(n)[2][i][1]: a_i:=ifactors(n)[2][i][2];
   x := x*(p_i^a_i-p_i^(a_i-1)+1): od: RETURN(x) end;
\end{verbatim}

\vskip1mm
\noindent {{\bf L\'aszl\'o T\'oth}\\
University of P\'ecs\\
Institute of Mathematics and Informatics\\
Ifj\'us\'ag u. 6\\
7624 P\'ecs, Hungary\\
ltoth@ttk.pte.hu}

\end{document}